\documentclass{article}
\usepackage{amsbsy,amssymb,amscd,amsfonts,latexsym,amstext,delarray,
amsmath, diagrams}  \headsep=15pt
\def\hpi{\hat{\pi}}

\def\cZ{ {\cal Z} }

\def\cX{{\cal X} }
\def\cD{ {\cal D}}

\def\bX{ \bar{X}}

\def\dim{{ \mbox{dim} }}
\def\Spec{{ \mbox{Spec} }}

\def\ra{{ \rightarrow }}

\def\g{{ \gamma }}

\def\hra{{ \hookrightarrow }}

\def\bs{ \backslash}

\def\G{{ \Gamma }}
\def\Gal{{ \mbox{Gal} }}

\def\bQ{\bar{\Q}}
\def\bF{{ \bar{F} }}
\def\bb{ \bar{b} }
\def\cE{ {\cal E}}

\def\Z{{ \mathbb{Z}}}

\def\bq{\begin{quote}}
\def\eq{\end{quote}}

\def\cJ{ {\cal J}}

\newtheorem{thm}{Theorem}[section]
\newtheorem{cor}[thm]{Corollary}
\newtheorem{lem}[thm]{Lemma}

\def\Q{\mathbb{Q}}

\def\cD{ {\cal D} }
\def\cY{ {\cal Y}}

\def\be{\begin{equation}}
\def\ee{\end{equation}}
\def\k{ \kappa}

\def\bY{\bar{Y}}

\def\bF{ \bar{F}}

\def\D{ \Delta}
\def\n{ {(n)}}

\def\s{ \sigma}
\def\bE{ \bar{E}}

\def\bF{ \bar{F}}

\def\bx{ \bar{x}}
\def\loc{\mbox{loc}}
\def\bb{\bar{b}}
\def\by{ \bar{y}}

\def\bZ{ \bar{Z}}
\def\cZ{ {\cal Z}}
\def\hP{ \hat{P}}
\title{The $l$-component of the unipotent Albanese map }
\author{Minhyong Kim and Akio Tamagawa}
\begin{document}
\maketitle
\begin{abstract}
We examine the local $l\neq p$-component of the $\Q_p$ unipotent Albanese map
for curves. As a consequence, we refine the non-abelian Selmer
varieties arising in the study of global points and deduce thereby a
new proof of Siegel's theorem for affine curves over $\Q$ of genus
one with Mordell-Weil rank at most one.
\end{abstract}
Let $F$ be a finite extension of $\Q_l$ for some
prime $l$, $R$ its ring of integers, and $k$ its residue
field. Let $\cX$ be a  curve over
$\Spec(R)$ with smooth generic fiber $X$.
 Let $\bb$ be an integral point,
i.e.  an $R$ section, of $\cX$. Given such a point, we will
denote by $b$ the induced $F$-point on $X$ and $b_s$
the reduction to a point of $\cX_s=\cX\otimes k$.
Let $\D=\hpi^{(l')}_1(\bX, b)$, the prime-to-$l$ geometric
fundamental group of $X$, where
we are viewing $b$ also as a geometric point coming in
from $\Spec(\bF)$. That is, $$\D=Aut(F_b)$$ and $F_b$ is
the fiber functor determined by $b$ on the category of finite \'etale covers
of $\bX$ whose Galois closures have degree prime to $l$.

Let $x\in X(F)$ be another  point  and denote by $\hP(x)$ the
$\D$ torsor of  prime-to-$l$ \'etale paths
$$\hP(x)=Isom(F_b,F_x)$$
in $\bX$ from $b$ to $x$. By choosing any path $p$ in $\hP(x)$ and measuring its
lack of $G_l$-invariance, that is, for each $g\in G_l$ we write
$g(p)=pc_g$  for a unique element $c_g\in \D$, we get a class $$[\hP(x)] \in H^1(G_l, \D)$$
$$g\mapsto c_g$$
 in a continuous non-abelian cohomology set for
 $G_l=\Gal(\bF/F)$. Recall that  non-abelian cohomology is
not a group, but has a well-defined base-point
$0\in H^1(G_l,\D)$ corresponding to the trivial torsor.
We can also carry out this construction for quotients of
$\D$ modulo its solvable series.
 That is, we define
$\D^{(1)}=\D$ and $\D^{(n+1)}=[\D^\n,\D^\n]$,
and then
$\D_\n=\D^\n\bs \D$. Each point $x$ determines also in an obvious way a
$\D_\n$-torsor $\hP_\n(x)$ which then ends up defining a class
$[\hP_\n(x)] \in H^1(G_l, \D_\n)$.
We wish to point out the following
\begin{thm}
Suppose $b_s$ is a smooth point of $\cX$ and
$\bx$ is an integral point such that $x_s=b_s$. Then
the class $[\hP_\n(x)] \in H^1(G_l, \D_\n)$ is zero.
\end{thm}

\begin{cor}
The map
$$\cX(R) \ra H^1(G_l, \D_\n)$$
$$x \mapsto [\hP_\n(x)]$$
has finite image.
\end{cor}

Our main motivation for writing down this easy result comes from  the structure theory of global non-abelian
Selmer varieties. That is, let $\cZ'$ be a proper
 curve over $\Spec(\Z)$ with good reduction outside a
finite set $S$ of primes. Let $\cD$ be a horizontal divisor
on $\cZ'$ and let $\cZ=\cZ'\setminus \cD$.
Fix a prime $p\notin S$ such that $\cD|(\cZ'\otimes \Z_{(p)})$
is \'etale over $\Spec (\Z_{(p)})$. Let $Z$ be the generic
fiber of $\cZ$.
Let $U^{et}=\pi_1^u(\bZ,b)$ be the pro-unipotent $\Q_p$-\'etale fundamental group
of $\bZ$ and $U^{et}_n=(U^{et})^n\bs U^{et}$ its quotient by the
descending central series, normalized so that
$(U^{et})^1=U^{et}$ and $(U^{et})^{n+1}=[U^{et},(U^{et})^n]$. We considered in
\cite{kim1}  and \cite{kim2} (section 2) the unipotent version $x \mapsto P^{et}_n(x)$
of the `torsor of paths' map that takes values in a
non-abelian Selmer variety
$$H^1_f(\G, U^{et}_n)$$
classifying torsors that are unramified outside $S$ and crystalline
at $p$. Here, $\G=\Gal (\bQ/\Q)$ and we are changing notation a bit
from the previous papers by including the condition of
being unramified outside $S$ into the subscript $f$. We thus have a sequence of maps
$$\k_n:\cZ(\Z) \ra H^1_f(\G, U^{et}_n)$$
that we called the {\em unipotent Albanese maps} (or alternatively, {\em unipotent Kummer maps})
which can be used to study the global integral points of $\cZ$.
According to the theorem and the corollary, the local components
of this map obtained by composing with the natural restriction maps
$$\loc_l:H^1_f(\G, U^{et}_n) \ra H^1(G_l, U^{et}_n)$$
to local Galois groups for various $l$ have finite image when $l\in S$.
It is easy to see that the image is zero for
$l\notin S$, $l\neq p$ (see section 2).

We specialize to the situation where we have a
proper model $\cE$  of an elliptic curve $E$ and $\cZ$ is obtained
by removing the closure of the origin.
In this case:
\begin{cor}
Assume $E(\Q)$ has rank $\leq 1$. Then
the Zariski closure of the image of
$$\k_3: \cZ(\Z) \ra H^1_f(\G, U^{et}_3)$$
has dimension $\leq 1$.
\end{cor}
We have stressed in the work cited the importance of conducting a refined study of
the unipotent Albanese map and its implications for the arithmetic of hyperbolic curves.
This discussion will not be repeated here, but
subsequent to the proof, we will remind the reader how this corollary implies
Siegel's theorem  for
$\cZ$ saying that $\cZ(\Z)$ is finite. We obtain thereby a new class of hyperbolic curves for which
the $\pi_1$-approach to finiteness yields definite results.
It is perhaps worthwhile to note that
the construction of the `refined' Selmer variety ($H^1_{\Sigma}$ in the next section) leading up to
the last corollary illustrates precisely the utility of the global
classifying space for torsors. The point is that there are
more natural algebraic conditions to impose on the global
classifying space than on the local one, giving us
control of the Zariski closure of interest.

\section{Proof of theorem}

Let $I$ be the directed category of natural numbers relatively prime
to $l$ where $i \ra k$ if $i|k$. It will be convenient to
deal first somewhat formally with the case $n=2$.
\begin{lem}
There is an inverse
system $Y_2(i)$ of finite \'etale covers of
$X$  over $F$ indexed by $I$ and equipped with:

(1)  smooth, $R$-models $\cY_2(i)$;

(2) a compatible system of \'etale maps
$$\begin{diagram}
 \cY_2(k) & &\rTo & & \cY_2(i) \\
 & \rdTo & & \ldTo \\
 & & \cX_{sm}& &
 \end{diagram} $$
for $i|k$, where
$\cX_{sm}$ is the smooth  locus of $\cX$; and

(3)  a compatible system of $R$ points
$$\bb_2(i): \Spec(R) \ra \cY_2(i)$$
lying over $\bb \in \cX(R)$;

such that the base change to $\bF$ is a cofinal
system for the category of abelian, prime-to-$l$ covers
of $\bX$.
\end{lem}
{\em Proof.} This is the standard generalized Jacobian construction
with a minor commentary about the models.
First, note that we can map
the base-point $b$ to the origin $e$ in the generalized Jacobian \cite{serre}
$J$ of $X$. The $Y_2(i)$ will then be
the pull-back to $X$ of the $i$-multiplication
on $J$. Now, let $\cJ$ be the N\'eron model of
$J$ (\cite{BLR}, theorem 10.2.2). By the universal property,
we have a map
$$\cX_{sm} \ra \cJ$$
from the smooth locus of $\cX$ to the N\'eron model.
The multiplication by $i$ clearly extends to an \'etale
 map of finite type (but not necessarily finite)
on $\cJ$. In any case, let $\cY_2(i)$ be
defined by the Cartesian diagram
$$\begin{diagram}
\cY_2(i) & \rTo & \cJ \\
\dTo & & \dTo_{i} \\
\cX_{sm} & \rTo & \cJ
\end{diagram}$$
which then fits into a compatible system. Fortunately,
the origin $\bar{e}:\Spec (R) \ra \cJ$ maps by $i$ to
$\bar{e}$ which is also  the image of $\bb$ in $\cJ$.
So we get a system of points
$$\bb_2(i)=(\bb, \bar{e}):\Spec (R) \ra \cY_2(i)$$
Since the maps between the $\cY_2(i)=\cX_{sm}\times_{\cJ}\cJ$ are induced by the
$i$-power maps on the second factor that preserve the origin,
these points are compatible.
$\Box$

Now, as an easy consequence of Hensel's lemma, we get:
\begin{cor}
The system $\cY_2(i)$ is equipped with a compatible collection of points
$$\bx_2(i):\Spec(R) \ra \cY_2(i)$$
lying over $\bx \in \cX(R)$ such that $x_2(i)_s=b_2(i)_s$.
\end{cor}

In particular, we have shown that the
class
$$[\hP_{(2)}(x)] \in H^1(G_l, \D_{(2)})$$
is trivial. We proceed now by induction in
an obvious way. That is,
we will show
that there exists a cofinal system
$\{ \bY_n(i)\}$  of Galois covers of $\bX$
for $\D_\n$
with $F$-models $$Y_n(i)\ra X,$$ smooth $R$
models
$\cY_n(i)$, and compatible \'etale diagrams
$$\begin{diagram}
\cY_n(k)& & \rTo & &\cY_n(i)\\
& \rdTo & & \ldTo & \\
& &  \cX & &
\end{diagram}$$
for $i|k$. There will be
two  compatible collections
$\{\bb_n(i)\}$ and $\{\bx_n(i)\}$ of $R$ points lying above
$\bb$ and $\bx$ such that $b_n(i)_s=x_n(i)_s$.  Notice
that once the maps exist they will automatically
be \'etale, since both source and target are asserted to be \'etale over $\cX$. Also,
if the base-point given by the $b_n(i)$ are compatible,
then the $x_n(i)$ are automatically compatible
since they lie over $x$ and reduce to the same
points as the $b_n(i)$.
We will
thereby prove the triviality
of $\hP_\n(x)$.

Assume that the
system
$\{Y_{n-1}(i)\}$ has already been constructed.
Denote by
$$Y_{n}(i)$$ the \'etale cover of $Y_{n-1}(i)$
obtained by pulling back from its Jacobian
$J_{n-1}(i)$  the isogeny
$$i:J_{n-1}(i) \ra J_{n-1}(i)$$
with the map that sends $b_{n-1}(i)$ to the origin.
The same argument using the N\'eron model as in the lemma
shows that we have $R$ models $\cY_{n}(i)$
\'etale over the
 $\cY_{n-1}(i)$, and hence, over $\cX$. Since $\bY_{n}(i)$ is defined
 by a characteristic subgroup of the prime-to-$l$ fundamental group of
 $\bY_{n-1}(i)$, it is Galois over $\bX$.
 Also, by
 construction, the base-point
 $\bb_{n}(i)$ maps to the basepoint $\bb_{n-1}(i)$,
 and hence, to $\bb\in \cX$. By Hensel's lemma we again
 get points
$\bx_{n}(i) \in \cY_n(i)$ that map to $\bx_{n-1}(i)$
and hence to $\bx$ and furthermore reduce to the same point
on the special fiber as $\bb_{n}(i)$. We need to check that these are compatible.
This follows from the existence of the $R$-integral transition maps
from the previous paragraph, given any two objects
$Y_{n}(i)$ and $Y_{n}(k)$ with $i|k$.
By induction, we already have the map from $\cY_{n-1}(k)$ to
$\cY_{n-1}(i)$ preserving basepoints.
Now we compare
$\cY_{n}(k)$ and $\cY_{n}(i)$.
Consider the diagram
$$
\begin{diagram}
\cY_{n}(k) & & \rDashto & &\cY_{n}(i) & &  \\
            & \rdTo & & &\vLine &\rdTo &     \\
\dTo & & \cJ_{n-1}(k) & \rTo^{(k/i)g} &\HonV & & \cJ_{n-1}(i) \\
     & &\dTo^{k} & & \dTo & & \\
\cY_{n-1}(k) & \hLine & \VonH &\rTo & \cY_{n-1}(i) & &\dTo^{i} \\
 &\rdTo & & & & \rdTo & \\
 & & \cJ_{n-1}(k) & & \rTo^g & & \cJ_{n-1}(i) \\
 \end{diagram} $$
 The right and left sides of the box are fiber diagrams
 and the map $g$ between the (Neron models of the) Jacobians is induced by the
 maps between the curves and the Neron model property,
 and hence, the bottom square commutes. The front
 square commutes because $g$ is
 a homomorphism. Therefore, by the universal property of
 the fiber product, the dotted arrow can be filled in so that
 the top and back squares commute. Finally, since the
 base point $b_{n}(k)$ maps to the
 origin in $J_{n-1}(k)$ and then from there to the the origin
 in $J_{n-1}(i)$ and to $b_{n-1}(k)$ which then
 maps to $b_{n-1}(i)$, the dotted arrow preserves base-points.
We are done with the construction of the cofinal system.

We remind ourselves briefly how this construction trivializes
$\hP_\n(x)$. Given any object $V \ra \bX$ in the category corresponding
to $\D_\n$, we define a bijection $$p_V:V_b\simeq V_x$$ as follows: Let $v\in V$.
Choose an object $(\bY_n(i), b_n(i))$  in our cofinal system with  a map
$$f_v:(\bY_n(i), b_n(i)) \ra (V,v)$$
Then by definition $p_V(v)=f_v(x_n(i))$. Suppose
we chose another such map
$$g_v:(\bY_n(j), b_n(j)) \ra (V,v)$$
Then we find $(\bY_n(k), b_n(k))$ as above with $ij|k$
that fits into a commutative diagram
$$\begin{diagram}
(\bY_n(k), b_n(k)) & \rTo^{f_v} & (\bY_n(i), b_n(i))\\
\dTo &\rdTo^{h_v} & \dTo \\
(\bY_n(j), b_n(j)) & \rTo^{g_v} & (V, v)
\end{diagram} $$
The reason the diagram has to commute is that $b_n(k)$
maps to $v$ both ways and $V$ is \'etale over $X$.
Thus, $p_V(v)=h_v(x_n(k))$ is independent of the choice. One checks
in a straightforward way that each $p_V$ is a bijection and that
it is compatible with maps between covering spaces.

We recall why  the path $p$ is then $G_l$-invariant,
that is
$$p_{\s^*(V)}( \s(v))=\s(p_V(v))$$
for any element $\s\in G_l$. Here,
we denote also by $\s$  the natural map of $F$-schemes
$$\s: V\ra \s^*(V)$$
inverse to the base-change map.
Let $$f_v(\bY_n(i),b_n(i)) \ra (V,v).$$
Then $$f_v^{\s}:(\bY_n(i),b_n(i)) \ra (\s^*V,\s(v))$$
where $f_v^{\s}=\s\circ f\circ \s^{-1}$ is now a map of
$\bF$-schemes. So then
$$p_{\s^*(V)}(\s(v))=f_v^{\s}(x)=\s(f_v(\s^{-1}(x)))=\s(f_v(x))$$
since $x$ is $G_l$-invariant. But $f_v(x)=p_V(v)$,
giving us the desired equality.
\medskip

Remark: It is possible to generalize the statement to the case where
$\D_\n$ is replaced by the entire $\D$ using the theory of log fundamental groups.
However, we will not present the proof here since the added generalization does
not enhance the application at the moment.

\section{Proof of corollaries}
{\em Corollary 0.2}

By resolution of singularities for two-dimensional schemes, we can assume $\cX$
is regular, and hence, that for any point $\bx\in \cX(R)$, its reduction $x_s$
is smooth. Then for any other point $\by$ such that $y_s=x_s$, we have
produced in the proposition a $G_l$-invariant path $\g$ in $\hpi_1^{(l)}(\bX; x,y)$. Composition
with $\g$ then induces an isomorphism of $\G_l$-equivariant torsors
$$\g\circ :\hpi_1^{(l)}(\bX;b,x)\simeq \hpi_1^{(l)}(\bX;b,y)$$
That is to say, the map
$$\cX(R) \ra H^1(G_l, \D_{\n})$$
factors through $$\cX(k),$$ the points in the residue field.
Hence, it has finite image. $\Box$

{\em Corollary 0.3}

We can again assume that $\cE$ is regular. Denote
by
$$\Sigma_n \subset \prod_{l\in S} H^1(G_l, U^{et}_n)$$
the image of  $\prod_{l\in S} \cZ(\Z_l)$.
This is finite by the previous corollary.

We recall briefly the situation at $l\notin S, l\neq p$. There, an argument
following the proof of the theorem
shows us that
the classes
$$\hP_{\n} \in H^1(G_l, \D_\n)$$
go to zero under the restriction map to the inertia group
$$H^1(G_l, \D_\n)\ra H^1(I_l, \D_\n)$$
This is because, using the N\'eron models with good reduction,
 the tower $\cY_n(i)$ can be taken proper and \'etale over $\cX$. Therefore,
we have isomorphisms of fiber functors
$$F_x \simeq F_{x_s}$$
induced by
the compatible bijections
$$Y_n(i)_x \simeq \cY_n(i)_{x_s}$$
showing that the $G_l$-action on all of the
$$Isom (F_b, F_x)$$
are unramified.
Therefore, the images of the classes $[P^{et}_n(x)]$ in
$$H^1(I_l, U^{et}_n)$$
are also trivial. Let us see
that this implies triviality in
$H^1(G_l, U^{et}_n)$.
We use the exact sequence
$$0\ra (U^{et})^{n+1}\bs (U^{et})^n \ra U^{et}_{n+1} \ra U^{et}_n \ra 0$$
By induction, the class
$P^{et}_{n+1}(x) \in H^1(G_l,U^{et}_{n+1} )$
comes from
$H^1(G_l,(U^{et})^{n+1}\bs (U^{et})^n)$.
But because it is unramified, it actually comes from
$$H^1(\Gal(\bar{k}/k), (U^{et})^{n+1}\bs (U^{et})^n)$$
On the other hand, there is a Galois equivariant surjection
$$\begin{diagram}
[(U^{et})^{2}\bs U^{et}]^{\otimes n} &\rOnto & (U^{et})^{n+1}\bs (U^{et})^n\end{diagram}$$
so an obvious weight argument shows that
$$H^1(\Gal(\bar{k}/k), (U^{et})^{n+1}\bs (U^{et})^n)=0$$
giving us the desired vanishing. One point of care in the discussion that we've presented
in this paper is that
the groups
$H^1(I_l,U^{et}_{n} )$
are not necessarily representable by varieties, at least using the technique of \cite{kim1}.
But because of the vanishing in $H^1(G_l,U^{et}_{n} )$, which are representable,
we do not need to deal with this
for $l$ of good reduction. (Alternatively and equivalently, we could have
started out with a restricted ramification variety $H^1(\G_T, U^{et}_n)$ as in
op. cit.)

To return to the proof of the corollary, for $m\leq n$,
define the subvariety
$$H^1_{\Sigma_m} (\G, U^{et}_n) \subset H^1_{f} (\G, U^{et}_n)$$
as the intersection of
$$(\prod_{l\notin S\cup \{p\}} \loc_l)^{-1}(0)$$
and
$$(p_{n,m}\circ \prod_{l\in S} \loc_l)^{-1}(\Sigma_m),$$
where $p_{n,m}$ denotes the projection
$$\prod_{l\in S} H^1(G_l, U^{et}_n) \ra \prod_{l\in S} H^1(G_l, U^{et}_m)$$
To see that the second condition does define a subvariety, one again proceeds as in loc. cit.
by first looking at $H^1_{\Sigma_2}(\G, U^{et}_n)$.
Since this is defined as an intersection of inverse images of
points under functorial maps to vector groups, it is a closed subvariety.
Now proceed by induction on $m$ and assume
$H^1_{\Sigma_{m-1}}(\G, U^{et}_n)$ is a subvariety.

For each of the finitely many $v \in \Sigma_{m-1}$ the fiber
$$[\prod H^1(G_l, U^{et}_m)]_v$$
is also represented by a product of vector spaces (with possibly different
dimensions for different $v$). And hence, the intersection
$$[\prod H^1(G_l, U^{et}_m)]_v\cap \Sigma_{m}$$
defines a subvariety (of points) and from there, we get the algebraicity
of
$$H^1_{\Sigma_m}(\G, U^{et}_n)$$

We have seen that
the Zariski closure of
$\k_3(\cZ(\Z))$ lies in
$H^1_{ \Sigma_3} (\G, U^{et}_3)$.
Now examine the sequence
$$0\ra H^1_f(\G, (U^{et})^{n+1}\bs (U^{et})^n) \ra H^1_f(\G,U^{et}_{n+1}) \ra H^1_f(\G,U^{et}_n)$$
which is exact in the sense that the vector group on the
left acts on the middle variety, with orbit space the image of the second map.
This will induce a sequence
$$0\ra H^1_{\Sigma_3}(\G, (U^{et})^{3}\bs (U^{et})^2) \ra H^1_{\Sigma_3}(\G,U^{et}_{3}) \ra H^1_{\Sigma_2}(\G,U^{et}_2)$$
which is exact in the naive sense that the inverse image of the base-point under the second
map is the image of the first map. But if we examine any other
fiber $H^1_{\Sigma_3}(\G,U^{et}_{3})_v$
of the second map for  $v\in H^1_{\Sigma_2}(\G,U^{et}_{2})$ as above, we see that it is contained in a set of the form
$$\tilde{v}+H^1_{\Sigma_{3,w}-\tilde{w}}(\G, (U^{et})^{3}\bs (U^{et})^2)$$
where $w=\prod_{\l\in S}\loc_l(v)$, $\Sigma_{3,w}$ is the set of points in $\Sigma_3$ mapping
to $ w$, and the tilde denotes liftings to $H^1_f(U^{et}_3)$.
Since $Z$ is an elliptic curve minus the origin,
we have
$$(U^{et})^{3}\bs (U^{et})^2\simeq \Lambda^2 ((U^{et})^2\bs U^{et}) \simeq \Q_p(1).$$
(The first isomorphism arises because there is a surjection from the
wedge product, but the two spaces have the same dimension (\cite{kim2}, section 3).)
Meanwhile, the classes in $H^1_f(\G, \Q_p(1))$
are already crystalline at $p$ and unramified outside of $S$,
forcing an injection (\cite{BK}, example 3.9)
$$H^1_f(\G, \Q_p(1))\hra \prod_{l\in S}H^1(G_l, \Q_p(1))$$
It is useful here to recall that this is an injection of varieties, not
just $\Q_l$ points, because both of the cohomologies are represented by
vector groups.
Therefore, each of the
$H^1_{\Sigma_{3, w}-\tilde{w}}(\G, \Q_p(1))$ are finite varieties.
We conclude that the
Zariski closure of
$\k_3(\cZ(\Z))$ is quasi-finite over
the closure of
$\k_2(\cZ(\Z))$. On the other hand, the latter is contained in
$\k_2(E(\Q))$ which is a subgroup of rank $\leq 1$ of the vector group
$H^1_f(\G, U^{et}_2)$. From this,
we get the desired dimension inequality
$$\dim \overline{\k_3(\Z)} \leq 1$$
$\Box$

To arrive at Siegel's theorem, we use  the map \cite{kim2}
$$H^1_f(\G, U^{et}_3) \stackrel{\loc_p}{\ra} H^1_f(G_p, U^{et}_3) \ra U^{dr}_3/F^0$$
The calculation of \cite{kim2}, section 4, gives
$\dim U^{dr}_3/F^0=2$
and hence, the image of
$\k_3(\Z)$ cannot be Zariski dense. Therefore, it must be finite as in \cite{kim1}, section 3,
following the non-abelian method of Chabauty.
\section{Comment}
In the abelian theory, say of the elliptic curve $E/\Q$, the {\em triviality}
of $$E(\Q_l) \ra H^1(G_l, H^{et}_1(\bE,\Q_p))$$
for $l\neq p$
plays an important role, even as the elementary nature of the
fact tends to obscure its significance. In fact, the proof in this case
is made even easier by the incompatibility of the $l$-adic and $p$-adic
group structures. Let us recall the application to finiteness:
if one has a global cohomology class
$c\in H^1(\G,H^{et}_1(\bE,\Q_p))$, then
the sum of local pairings gives us
$$\Sigma_l<\loc_l(c), \loc_l(\k_2(P))>_l=0$$
for every global point $P$ by local-global duality. On the other hand,
because $\loc_l(\k_2(P))=0$ for all $l\neq p$, this becomes
the single equation
$$<\loc_p(c), \loc_p(\k_2(P))>_p=0$$
implying the finiteness of $E(\Q)$ whenever
$\loc_p(c)$ is not in the finite part $H^1_f(G_p, H^{et}_1(\bE,\Q_p))$.

In our non-abelian situation, because of the lack of a group
structure, we only get finiteness of the
image of $\loc_l\circ \k_n$ for $l\neq p, l\in S$ as soon as $n\geq 3$.
It is something of an interesting problem to figure out
if the image can be made trivial, for example, if
$p$ is taken sufficiently large. (It {\em can} be made trivial on the
points that reduce to the same connected component as the basepoint
in the special fiber of
a regular model, although we
will not dwell here on this fact.)
But if one imagines some sort of a `non-abelian duality'
to play a role in finiteness, the result we have should
already be sufficient. What we have in mind here
is the existence of a global object $c$ together with
localizations $\loc_l(c)$ each of which are algebraic functions
on $H^1(G_l, U^{et}_n)$.
The analogue of local-global duality should then provide a
statement of the form
$$\Sigma_l\loc_l(c)[\loc_l(\k_n(P))]=0$$
for every global point $P$.
But using the finite set $\Sigma$,
the elements $$\s=(\s_l) \in \Sigma\subset \prod_{l\in S} H^1(G_l, U^{et}_n)$$
 give us a finite collection of equations
$$\loc_p(c)[\loc_p(\k_n(P))]=-\Sigma_{l\in S}\loc_l(c)[\s_l]$$
at least one of which has to be satisfied by
every global point. That is to say, provided such a conjectural
framework can actually be realized, we need only replace a single
equation for $\loc_p(\k_n(P))$ by {\em finitely many }
equations.
\medskip

{\bf Acknowledgements:}
\smallskip

M.K. was supported in part by a grant from the National Science Foundation and
a visiting professorship at RIMS.
He is  extremely grateful to Kazuya Kato, Shinichi Mochizuki, A. T.,
and the staff at RIMS for providing the stimulating environment
in which this work was completed.

{\footnotesize M.K.: Department of Mathematics, Purdue University, West Lafayette, IN 47906, U.S.A. and
Department of Mathematics, University of Arizona, Tucson, AZ 85721, U.S.A. e-mail: kimm@math.purdue.edu}

{\footnotesize A.T.: Research Institute for Mathematical Sciences, Kyoto University, Kyoto 606-8502, Japan.
e-mail: tamagawa@kurims.kyoto-u.ac.jp}

\end{document}